\newcommand{\bull}{\vrule height .9ex width .8ex depth -.1ex}
\newcommand{\ppp}{\hfill $\bull$ }
 \documentclass[11pt]{article}
 \author{ M.-L. Labbi\thanks{
  Address: Department of Mathematics, College
 of Science, University of Bahrain, Isa Town 32038 Bahrain.
  E-mail: labbi@sci.uob.bh }}
   \title{Double forms, curvature structures and the (p,q)-curvatures}
   \date{}
   \usepackage{amsmath}
   \usepackage{amscd}
\newtheorem{theorem}{Theorem}[section]
\newtheorem{corollary}[theorem]{Corollary}
\newtheorem{lemma}[theorem]{Lemma}

\newtheorem{remark}{Remark}[section]
\newtheorem{proposition}[theorem]{Proposition}
\newtheorem{definition}{Definition}[section]

\begin{document}
   \maketitle
   \begin{abstract} We introduce a  natural extension of the
   metric tensor and the Hodge star operator to the algebra of double forms to study some aspects 
of the structure of this algebra.  These
   properties are then used to
    study new Riemannian curvature invariants, called the $(p,q)$-curvatures. They  are 
a generalization of the $p$-curvature obtained by substituting the
Gauss-Kronecker tensor to the Riemann curvature tensor. In particular, for $p=0$, the $(0,q)$-curvatures coincide with the H. Weyl curvature invariants, for $p=1$ the $(1,q)$-curvatures are the curvatures of generalized Einstein tensors and for $q=1$ the $(p,1)$-curvatures coincide with the p-curvatures.\\
Also,  we
prove that for an Einstein manifold of dimension $n\geq 4$ the
second H. Weyl curvature invariant is nonegative, and that it is
nonpositive for a conformally flat manifold with zero scalar
curvature. A similar result is proved for the higher H. Weyl curvature invariants.

   \end{abstract}
 \section{ Introduction}
Let $(M,g)$ be a smooth Riemannian manifold of dimension n. We denote by
 $\Lambda^{*}M=\bigoplus_{p\geq 0}\Lambda^{*p}M$ the ring of differential
  forms
 on $M$.  Considering the tensor product over the ring of  smooth functions,
  we define
 ${\cal D}= \Lambda^{*}M\otimes \Lambda^{*}M=\bigoplus_{p,q\geq 0}
  {\cal D}^{p,q}$ where $  {\cal D}^{p,q}= \Lambda^{*p}M \otimes
   \Lambda^{*q}M$.
  It is  graded associative  ring and  called the ring of
   double forms on $M$. \par
   The ring of curvature structures on $M$ is the ring ${\cal C}=\sum_{p\geq 0}
   {\cal C}^p$ where ${\cal C}^p$ denotes symmetric elements in
   ${\cal D}^{p,p}$.\par
   These notions have been developed  by Kulkarni \cite{Kulk},
   Thorpe \cite{Thorpe} and other mathematicians.\par\noindent
   The object of this paper is to  study   some  properties of
   these structures, and then to apply them to study generalized $p$-curvature
   functions.\par\medskip
    The paper is divided into 5 sections. In section 2, we study the
    multiplication map  by $ g^l$   in ${\cal D}^{p,q}$. In particular
    we prove that it is one to one if $p+q+l\leq n$. This result will
    play an important role in simplifying complicated calculations (see
    for example the applications in sections 5).
    We also  deduce some properties of the multiplication map by $g$.
\par\medskip
In section 3, we introduce a natural inner product $<,>$ in ${\cal D}$ and
we extend the Hodge operator $*$ in a natural way to ${\cal D}$. Then, we prove
two simple relations between the contraction map and the
 multiplication map by $g$, namely for all $\omega \in {\cal D}$,we have
 $$g\omega=*c*\omega$$
 Also we prove that the contraction map is the adjoint of the multiplication by
 $g$, precisely  for all $\omega_1, \omega_2\in {\cal D}$, we have
 $$<g\omega_1,\omega_2>=<\omega_1,c\omega_2>$$
 and we deduce some properties of the contraction map c.\par
 At the end of this section, we deduce a canonical orthogonal decomposition
 of ${\cal D}^{p,q}$ and we give explicit formulas for the
 orthogonal projections onto the different factors.\par\medskip
 In section 4, we concentrate on the ring
  of symmetric double forms satisfying the first
 Bianchi identity which shall be denoted  by ${\cal C}_1$. We prove in this context
 a useful explicit
 formula for the  Hodge star operator. Also, we emphasize its action on
 the different factors of the previous orthogonal decomposition of
 double forms.\par\medskip
In section 5, we define new Riemannian curvature invariants, namely the $(p,q)$- curvature tensors 
$R_{(p,q)}$ and their sectional
curvatures $s_{(p,q)}$.\\
Note that these curvatures include  many of well known curvatures.\\
For $q=1$, the $(p,1)$-curvature coincide with
 the $p$-curvature. In particular,  $s_{(0,1)}$ is the half of
  the scalar curvature and  $s_{(n-2,1)}$ is
   the sectional curvature of $(M,g)$.\\
For $p=0$ and $2q=n$, $s_{(0,{n\over 2})}$ is up to a constant
 the Killing-Lipshitz curvature. More generally, $s_{(n-2q,q)}(P)$ is,
  up to a constant, the Killing-Lipshitz curvature
of $P^{\bot}$.\par
For $p=0$, $s_{(0,q)}$
 are scalar functions which generalize the usual scalar curvature.
 They are, up to a constant, the H. Weyl curvature invariants, 
 that is the integrands in the Weyl tube formula \cite{Wey}.
\par\noindent
Finally, for $p=1$, $R_{(1,q)}$ are 
 generalized Einstein tensors. 
In particular, for $q=1$ we recover the usual Einstein tensor. \par\smallskip
This section contains also several examples and properties of these curvature invariants.
In particular we prove using the $(p,1)$-curvatures a characterization
 of Einstein metrics   and conformally flat metrics  with constant scalar
 curvature. Also a generalization of the previous result to the higher
  $(p,q)$-curvatures is proved.\par\medskip
  In the last section we prove under certain geometric hypothesis on the
   metric,
  a restriction on the sign of the  H. Weyl curvature invariants,
  that are the integrands
  in his well known tube formula \cite{Wey}.
    In particular we
   prove the following results:\par\medskip
{\sl If $(M,g)$ is an Einstein manifold with dimension $n\geq 4$,
then $h_4\geq 0$ and $h_4\equiv 0$ if and only if $(M,g)$ is flat.}
\par\medskip
{\sl If $(M,g)$ is a conformally flat  manifold with zero scalar curvature
and
 dimension $n\geq 4$,
then $h_4\leq 0$ and $h_4\equiv 0$ if and only if $(M,g)$ is flat.}

where $h_4$ is the  second H. Weyl curvature invariant , which can be defined by
$$h_4=|R|^2-|c(R)|^2+{1\over 4}|c^2(R)|^2$$
\par\noindent
where $R$ denotes the Riemann curvature tensor of $(M,g)$.

\section{The Algebra of double forms}

Let $(V,g)$ be an Euclidean real vector space  of dimension n. In the
 following
we shall identify whenever convenient (via their Euclidean structures),
the vector spaces
 with their duals. Let
  $\Lambda^{*}V=\bigoplus_{p\geq 0}\Lambda^{*p}V$ (resp.
   $\Lambda V=\bigoplus_{p\geq 0}\Lambda^{p}V$) denote the exterior algebra
 of $p$-forms (resp. $p$-vectors) on $V$. Considering the tensor product,
  we define the space of double forms
 ${\cal D}= \Lambda^{*}V\otimes \Lambda^{*}V=\bigoplus_{p,q\geq 0}
  {\cal D}^{p,q}$ where $  {\cal D}^{p,q}= \Lambda^{*p}V \otimes
   \Lambda^{*q}V$.  It is a bi-graded associative  algebra,
    where the multiplication is denoted by  a dot., we shall omit it whenever
     it is possible.\par\noindent
    For $\omega_1=\theta_1\otimes \theta_2\in { \cal D}^{p,q}$ and
    $\omega_2=\theta_3\otimes \theta_4\in  {\cal D}^{r,s}$, we have
    \begin{equation}
    \label{def:prod}
     \omega_1.\omega_2= (\theta_1\otimes \theta_2 ).(\theta_3\otimes
     \theta_4)=
    (\theta_1\wedge \theta_3 )\otimes(\theta_2\wedge \theta_4)\in
    {\cal D}^{p+r,q+s}\end{equation}
    Recall that each element of the tensor product
    $  {\cal D}^{p,q}= \Lambda^{*p}V \otimes \Lambda^{*q}V$
     can be identified canonically
     with a bilinear form $\Lambda^pV\times\Lambda^qV\rightarrow {\bf R}$. That is
     a multilinear form which is skew symmetric in the first $p$-arguments and also
     in the last $q$-arguments.
      Under this
      identification, we have

\begin{equation}
\begin{split}
\label{eps:prod}
 &\omega_1.\omega_2(x_1\wedge...\wedge x_{p+r},y_1\wedge...\wedge y_{q+s})\\
&= (\theta_1\wedge \theta_3 )(x_1\wedge...\wedge x_{p+r})
(\theta_2\wedge \theta_4)(y_1\wedge...\wedge y_{q+s})\\
&= {1\over p!r!s!q!}\sum_{\sigma\in S_{p+r}, \rho\in S_{q+s}}
\epsilon(\sigma)\epsilon(\rho)
\omega_1(x_{\sigma(1)}\wedge...\wedge x_{\sigma(p)};y_{\rho(1)}
\wedge...\wedge y_{\rho(q)})\\
&\phantom{...mmmmmmmmmmm}
\omega_2(x_{\sigma(p+1)}\wedge...\wedge x_{\sigma(p+r)};y_{\rho(q+1)}
\wedge...\wedge y_{\rho(q+s)})
\end{split}
\end{equation}
A similar calculation shows that
\begin{equation}
\begin{split}
\label{eps:power}
 \omega_1^k&(x_1\wedge...\wedge x_{kp},y_1\wedge...\wedge y_{kq})=\\
 &=(\theta_1\wedge ...\wedge\theta_1)(x_1\wedge...\wedge x_{kp})
  (\theta_2\wedge ...\wedge\theta_2)(y_1\wedge...\wedge y_{kp}) \\
 &={1\over (p!)^k(q!)^k}\sum_{ \sigma\in S_{kp}, \rho\in S_{kq}}
\epsilon(\sigma)\epsilon(\rho)
 \omega_1(x_{\sigma(1)}\wedge...\wedge x_{\sigma(p)};y_{\rho(1)}
\wedge...\wedge y_{\rho(q)})\\
 &\phantom{..........} ...\omega_1(x_{\sigma(p(k-1)+1)}\wedge...
\wedge x_{\sigma(kp)};y_{\rho(q(k-1)+1)}
\wedge ...\wedge y_{\rho(kq)})
\end{split}
\end{equation}
In particular, if $\omega_1\in {\cal D}^{1,1}$ we have
\begin{equation}
\label{h:power}
 \omega_1^k(x_1\wedge...\wedge x_{k},y_1\wedge...\wedge y_{k})=k!\det
 [\omega_1(x_i,y_j)]\end{equation}
We now introduce a basic map on ${\cal D}$:
\begin{definition} The  contraction $c$ maps  ${\cal D}^{p,q}$ into  ${\cal D}^{p-1,q-1}$.
Let $\omega \in  {\cal D}^{p,q}$, set $c\omega=0$ if $p=0$ or $q=0$.
Otherwise  set
$$c\omega(x_1\wedge...\wedge x_{p-1},y_1\wedge...\wedge y_{q-1})=
\sum_{j=1}^{n}\omega(e_j\wedge x_1\wedge... x_{p-1},
e_j\wedge y_1\wedge...\wedge y_{q-1})$$
where $\left\{e_1,...,e_n\right\}$ is an  orthonormal basis of $V$.
\end{definition}
The contraction map $c$ together with the multiplication map by $g$, (which
 shall be denoted also  by $g$), play a very important
  role in our study.
  \par\medskip\noindent
Let $\omega\in {\cal D}^{p,q}$, the following formula was proved
 in \cite{Kulk}
 \begin{equation}\label{kulk:form}
 c(g\omega)=gc\omega +(n-p-q)\omega
 \end{equation}
After consecutive applications of the previous formula, we get
\begin{equation*}
\begin{split}
c^k(g\omega)&=gc^k\omega +k(m+k-1)c^{k-1}\omega ,\, {\rm where}\, m=n-p-q\\
c^k(g^2\omega)&=g^2c^k\omega +2k(m+k-2)gc^{k-1}\omega+k(k-1)(m+k-3)(m+k-2)c^{k-2}\omega\\
c^k(g^3\omega)&=g^3c^k\omega +3k(m+k-3)g^2c^{k-1}\omega+3k(k-1)
(m+k-3)(m+k-4)gc^{k-2}\omega\\
&\phantom{.............}+k(k-1)(k-2)(m+k-3)(m+k-4)(m+k-5)c^{k-3}\omega\\
\end{split}
\end{equation*}
More generally, we have
\begin{lemma} For all $k,l\geq 1$ and $\omega \in D^{p,q}$,  we have
\begin{equation}
\label{ck:gl}
c^k({g^l\over l!}\omega)={g^l\over l!}c^k\omega
+\sum_{r=1}^{\min\{k,l\}}C_r^k \prod_{i=0}^{r-1}(n-p-q+k-l-i)
{g^{l-r}\over (l-r)!}.c^{k-r}\omega
\end{equation}
\end{lemma}
\begin{corollary} \label{ckgk:gkck}
If $n=p+q$ and $\omega\in D^{p,q}$ then
 for all $k$
we have
$$c^k(g^k\omega)=g^k(c^k\omega)$$
\end{corollary}
{\bf Proof.} After taking $k=l$ and $n=p+q$ in formula \ref{ck:gl}, we get
$$c^k({g^k\over k!}\omega)={g^k\over k!}.c^k\omega
+\sum_{r=1}^k C_r^k \prod_{i=0}^{r-1}(-i){g^{k-r}\over (k-r)!}c^{k-r}\omega
={g^k\over k!}c^k\omega$$
\ppp
\par\medskip\noindent
As a second consequence of the previous lemma, we get the following result
which  generalizes another lemma of Kulkarni \cite{Kulk}:\par\noindent
\begin{proposition}\label{mult:injec} The multiplication by $g^l$ is injective on $D^{p,q}$
whenever $p+q+l<n+1$.\end{proposition}
{\bf Proof.}\,  This property is true for $l=0$, suppose that
$g^{l-1}\omega =0\Rightarrow \omega =0$ for $p+q+l-1<n+1$. Suppose $g^{l}\omega=0$ and  $p+q+l<n+1$,  then the  contractions $c^k(g^l.\omega)=o$ for all $k$.\par
Taking $k=1,2,...,k,...,{\min\{p,q\}},{\min\{p,q\}}+1$
 and after a simplification (if needed) by
  $g^{l-1},g^{l-2},...,g^{l-k},...,g^{l-{\min\{p,q\}}+1}$,  respectively, we get using the previous lemma
\begin{equation*}
\begin{split}
-gc\omega &=l(n-p-q+1-l)\omega\\
-gc^2\omega&=(l+1)(n-p-q+2-l)c\omega\\
&\vdots\\
-g.c^k\omega &=(l+k-1)(n-p-q+k-l)c^{k-1}\omega\\
&\vdots\\
-g.c^{\min\{p,q\}}\omega &=(l+{\min\{p,q\}}-1)(n-{\max\{p,q\}}-l)
c^{\min\{p,q\}-1}\omega\\
0&=(l+{\min\{p,q\}})(n-{\max\{p,q\}}+1-l)c^{\min\{p,q\}}\omega
\end{split}
\end{equation*}
Consequently, we have $c^{\min\{p,q\}}\omega=...c^k\omega...=\omega=0$.
\ppp\par\noindent
\begin{remark}\label{rem:glck}\begin{enumerate}
\item The previous proposition cannot be obtained directly from
Kulkarni's Lemma\cite{Kulk}. Since consecutive applications of that lemma show that
the multiplication by $g^l$ is 1-1 only if $p+q+2l-2<n$.\par
\item We deduce from the previous proof that more generally we have
$$g^l\omega=0\Rightarrow c^k\omega=0\qquad {\mathrm for }\quad l+p+q< n+1+k$$
\end{enumerate}
\end{remark}
\begin{corollary}
\label{isom:cor}\begin{enumerate}
\item Let  $p+q=n-1$, then for each $i\geq 0$,
 the multiplication map by $g^{2i+1}$
$$g^{2i+1}:D^{p-i,q-i}\rightarrow D^{p+i+1,q+i+1}$$
is an isomorphism. In particular, we have
$$D^{p+i+1,q+i+1}= g^{2i+1}D^{p-i,q-i}$$
\item Let  $p+q=n$, then for each $i\geq 0$,
 the multiplication map by $g^{2i}$
$$g^{2i}:D^{p-i,q-i}\rightarrow D^{p+i,q+i}$$
is an isomorphism. In particular, we have
$$D^{p+i,q+i}= g^{2i}D^{p-i,q-i}$$
\end{enumerate}
\end{corollary}
{\bf Proof.} From the previous proposition it is 1-1. And for dimensions
reason (since $C_{p-i}^nC_{q-i}^n=C_{p+i+1}^nC_{q+i+1}^n$ if $p+q=n-1$, and
$C_{p-i}^nC_{q-i}^n=C_{p+i}^nC_{q+i}^n$ if $p+q=n$) then it is an
isomorphism.\ppp\par\medskip\noindent
The following proposition gives more detail about the multiplication by $g$
\begin{proposition}\label{gonto:bijective} The multiplication map by $g$ on $D^{p,q}$ is\par
\begin{enumerate}
\item one to one  if and only if $p+q\leq n-1$.
\item bijective if and only if $p+q=n-1$.
\item onto if and only if $p+q\geq n-1$.
\end{enumerate}
\end{proposition}
{\bf Proof.}
The only if part of the propostion is due simply to dimension
 reasons, so that parts 1) and 2) are direct consequences of Kulkarni's Lemma.
 \par
 Let now $i\geq 0$, $p_0+q_0=n-1$ for some $p_0,q_0\geq 0$ and
 $$g:D^{p_0+i,q_0+i}\rightarrow D^{p_0+i+1,q_0+i+1}$$
Remark that the restriction of the map $g$ to the subspace
$g^{2i}D^{p_0-i,q_0-i}$ of  $D^{p_0+i,q_0+i}$
is onto since its image is exactly   $g^{2i+1}D^{p_0-i,q_0-i}=D^{p_0+i+1,q_0+i+1}$
 by the previous proposition. The proof is similar in case there exists
$p_0,q_0\geq 0$ such that $p_0+q_0=n$. This completes the proof of
 the proposition.\ppp
 \par\medskip\noindent
\section{The natural Inner product and the Hodge star operator on $D^{p,q}$}
\subsection{The natural Inner product on $D^{p,q}$}
The natural metric on $\Lambda^{*p}V$
  induces in a standard way an inner product on
  $D^{p,q}=\Lambda^{*p}V \otimes \Lambda^{*q}V$. We shall denote it by $ <,>$.
  \par
We extend $<,>$ to ${\cal D}$ by declaring that  $D^{p,q} \perp D^{r,s}$
if $p\not = r$ or if $q\not =s$.
\begin{theorem}\label{theo:gc} If $\omega_1, \omega_2\in {\cal D}$  then
\begin{equation}
\label{adj:gc}
<g\omega_1,\omega_2>=<\omega_1,c\omega_2>
\end{equation}
that is the adjoint of the multiplication by $g$ is the  contraction map $c$.
\end{theorem}
{\bf Proof.} Let $\{e_1,...,e_n\}$ be an orthonormal basis of $V^*$. Since
the contraction map $c$ and the multiplication by g are
 linear it suffices to prove the theorem for
 $$\omega_2 =e_{i_1}\wedge ...\wedge e_{i_{p+1}}\otimes
  e_{j_1}\wedge ... \wedge e_{j_{q+1}} {\text{  and}} \,
 \omega_1 =e_{k_1}\wedge ...\wedge e_{k_{p}}\otimes
  e_{l_1}\wedge ...\wedge e_{l_{q}}$$
  where $i_1<...<i_{p+1};j_1<...<j_{q+1}; k_1<...<k_{p}$ and $l_1<...<l_{q}$.
Then
$$g\omega_1=\sum_{i=1}^n e_i\wedge e_{k_1}\wedge ...\wedge e_{k_{p}}
\otimes e_i\wedge e_{l_1}\wedge ...\wedge e_{l_{q}}$$
therefore
$$<g\omega_1,\omega_2>= \sum_{i=1}^n< e_i\wedge e_{k_1}\wedge ...
\wedge e_{k_{p}},e_{i_1}\wedge ...\wedge e_{i_{p+1}}>
<e_i\wedge e_{l_1}\wedge ...\wedge e_{l_{q}},
e_{j_1}\wedge ... \wedge e_{j_{q+1}}>$$
So it is zero unless if
$$e_{k_1}\wedge ... \wedge e_{k_{p}}=e_{i_1}\wedge ...
\hat{e}_{i_r}...\wedge e_{i_{p+1}}$$
$$e_{l_1}\wedge ... \wedge e_{l_{q}}=e_{j_1}\wedge ...
\hat{e}_{j_s}...\wedge e_{j_{q+1}}$$
and $i_r=j_s$ for some $r,s$. So that in this case we have
\begin{equation*}
\begin{split}
<g\omega_1,\omega_2>&= \sum_{i=1}^n< e_i\wedge e_{i_1}\wedge ...
\hat{e}_{i_r}...
\wedge e_{i_{p+1}},e_{i_1}\wedge ...\wedge e_{i_{p+1}}>\\
&\phantom{mmmmmmm}<e_i\wedge e_{j_1}\wedge ...\hat{e}_{j_s}..\wedge e_{j_{q+1}},
e_{j_1}\wedge ... \wedge e_{j_{q+1}}>\\
&=(-1)^{r+s}
\end{split}\end{equation*}
On the other hand, we have
$$c\omega_2 =0 \,  {\text{  if}} \,
\left\{i_1,...,i_{p+1}\right\}\cap \left\{j_1,...,i_{q+1}\right\}=\emptyset$$
otherwise,
$$c{\omega}_2 =
\sum_{{\scriptstyle i_r=j_s\atop\scriptstyle 1\leq r\leq p+1}\atop\scriptstyle 1\leq s\leq q+1}
(-1)^{r+s}
e_{1_1}\wedge ...\hat{e}_{i_r}..\wedge e_{i_{p+1}}\otimes
e_{j_1}\wedge ...\hat{e}_{j_s}..\wedge e_{j_{q+1}}$$
and therefore
\begin{equation*}
\begin{split}
<\omega_1,c\omega_2>=&
\sum_{{\scriptstyle i_r=j_s\atop\scriptstyle 1\leq r\leq p+1}\atop\scriptstyle 1\leq s\leq q+1}
(-1)^{r+s}
<e_{k_1}\wedge ...\wedge e_{k_{p}},
<e_{i_1}\wedge ...\hat{e}_{i_r}..\wedge e_{i_{p+1}}>\\
&\phantom{mmmmmmmmmm}< e_{l_1}\wedge ...\wedge e_{l_{q}}\otimes
e_{j_1}\wedge ...\hat{e}_{j_s}..\wedge e_{j_{q+1}}
\end{split}\end{equation*}
which is zero unless if
$$e_{k_1}\wedge ... \wedge e_{k_{p}}=e_{i_1}\wedge ...
\hat{e}_{i_r}...\wedge e_{i_{p+1}}$$
$$e_{l_1}\wedge ... \wedge e_{l_{q}}=e_{j_1}\wedge ...
\hat{e}_{j_s}...\wedge e_{j_{q+1}}$$
and $i_r=j_s$ for some $r,s$. In such case it is $(-1)^{r+s}$.
This completes the proof.\ppp

\subsection{Hodge star operator}
The Hodge star operator $*:\Lambda^{p}V^*\rightarrow \Lambda^{n-p}V*$
 extends in a natural way to a linear operator $*:{\cal D}^{p,q}
 \rightarrow {\cal D}^{n-p,n-q}$. If $\omega=\theta_1\otimes
\theta_2$ then we define $$*\omega=*\theta_1\otimes *\theta_2$$
Note that $*\omega(.,.)=\omega(*.,*.)$ as a bilinear form. Many properties
of the ordinary Hodge star operator can be extended to this new operator.
 We prove some of them below:
\begin{proposition} For all  $\omega,\theta\in D^{p,q}$, we have
\begin{equation}
\label{in:star}
<\omega,\theta>=*(\omega.*\theta)=*(*\omega.\theta)
\end{equation}
\end{proposition}
{\bf Proof.} Let $\omega=\omega_1\otimes\omega_2$ and $\theta=\theta_1\otimes
\theta_2$, then
\begin{equation*}
\begin{split}
\omega.*\theta&=(\omega_1\wedge *\theta_1)\otimes(\omega_2\wedge *\theta_2)\\
&=<\omega_1,\theta_1>*1\otimes <\omega_2,\theta_2>*1\\
&=<\omega,\theta>*1\otimes *1
\end{split}
\end{equation*}
This completes the proof. \ppp
\par\smallskip\noindent
The proof of the following properties is similar and straightforward
\begin{proposition}\begin{enumerate}
\item For all $p,q$, on $D^{p,q}$ we have
$$**=(-1)^{(p+q)(n-p-q)}Id$$
where $Id$ is the identity map on $D^{p,q}$.
\item For all $\omega_1\in D^{p,q},\omega_2\in D^{n-p,n-q}$ we have
$$<\omega_1,*\omega_2>=(-1)^{(p+q)(n-p-q)}<*\omega_1,\omega_2>$$
\item If $\bar \omega:\Lambda^p\rightarrow \Lambda^p$
 denotes the linear operator corresponding to $\omega\in D^{p,p}$, then
  $$*\bar\omega *:\Lambda^{n-p}\rightarrow \Lambda^{n-p}$$
 is the linear operator corresponding to $*\omega\in D^{n-p,n-p}$.
\end{enumerate}
\end{proposition}
 Using the Hodge star operator we can provide a nice formula relating the
 multiplication by $g$ and the contraction map $c$, as follows:
\begin{theorem} For every $\omega \in  D^{p,q}$, we have
\begin{equation}
\label{gstar:c}
g\omega=*c*\omega
\end{equation}
that is the following diagram is commutative for all $p,q$
\begin{equation*}
\begin{CD}
D^{p,q}@>g>>D^{p+1,q+1}\\
@VV{*}V     @AA{*}A\\
D^{n-p,n-q}@>c>>D^{n-p-1,n-q-1}
\end{CD}
\end{equation*}
\end{theorem}
\noindent
{\bf Proof.} The proof is similar to the one of theorem \ref{theo:gc}. Let
$\{e_1,...,e_n\}$ be an orthonormal basis of $V^*$, and let
 $$\omega =e_{i_1}\wedge ...\wedge e_{i_{p}}\otimes
  e_{j_1}\wedge ... \wedge e_{j_{q}}$$
  then
$$g\omega=\sum_{i=1}^n e_i\wedge e_{i_1}\wedge ...\wedge e_{i_{p}}
\otimes e_i\wedge e_{j_1}\wedge ...\wedge e_{j_{q}}$$
On the other hand, we have
$$*\omega=\epsilon(\rho)\epsilon(\sigma)e_{i_{p+1}}\wedge ...\wedge e_{i_{n}}
\otimes e_{j_{q+1}}\wedge ...\wedge e_{j_{n}}$$
so that
$$c*\omega =
\sum_{{\scriptstyle i_r=j_s\atop\scriptstyle p+1\leq r\leq n}\atop\scriptstyle q+1\leq s\leq n}
(-1)^{r+s}\epsilon(\rho)\epsilon(\sigma)
e_{i_{p+1}}\wedge ...\hat{e}_{i_r}..\wedge e_{i_n}\otimes
e_{j_{q+1}}\wedge ...\hat{e}_{j_s}..\wedge e_{j_{n}}$$
Therefore
\begin{equation*}
\begin{split}
*c*\omega & =
\sum_{{\scriptstyle i_r=j_s\atop\scriptstyle p+1\leq r\leq n}\atop\scriptstyle q+1\leq s\leq n}
(-1)^{r+s}\epsilon(\rho)\epsilon(\sigma)
*e_{i_{p+1}}\wedge ...\hat{e}_{i_r}..\wedge e_{i_n}\otimes
*e_{j_{q+1}}\wedge ...\hat{e}_{j_s}..\wedge e_{j_{n}}\\
&=
\sum_{{\scriptstyle i_r=j_s\atop\scriptstyle p+1\leq r\leq n}\atop\scriptstyle q+1\leq s\leq n}
{e}_{i_r}\wedge e_{i_{1}}\wedge ...\wedge e_{i_p}\otimes
{e}_{j_s}\wedge e_{j_{1}}\wedge ...\wedge e_{j_{q}}\\
&=g\omega
\end{split}\end{equation*}
This completes the proof.\ppp\par\medskip\noindent
As a direct consequence of the previous theorem and proposition
\ref{gonto:bijective}, we have
\par\smallskip\noindent
\begin{corollary} The contraction map $c$ on $D^{p,q}$ is
\begin{enumerate}
\item onto  if and only if $p+q\leq n-1$.
\item bijective if and only if $p+q=n-1$.
\item one to one if and only if $p+q\geq n-1$.
\end{enumerate}
\end{corollary}
\begin{corollary}\label{orth:decomp} For all $p,q\geq 0$ such that $p+q\leq n-1$,
 we have the orthogonal decomposition
$$D^{p+1,q+1}={\rm Ker}\, c\oplus gD^{p,q}$$
where $c:D^{p+1,q+1}\rightarrow D^{p,q}$ is the contraction map.
\end{corollary}
{\bf Proof.} First note that if $\omega_1\in \ker c$ and $g\omega_2\in g
D^{p,q}$ then by formula \ref{adj:gc}, we have
$$<\omega_1,g\omega_2>=<c\omega_1,\omega_2>=0$$
Next since $g$ is one to one and $c$ is onto, we have
$$\dim (gD^{p,q})=\dim D^{p,q}=\dim ({\rm image}\, c)$$
This completes the proof.\ppp
\begin{remark} \begin{enumerate}
 \item If $p+q>n-1$ , then  we have $\ker c=0$
and $D^{p+1,q+1}$  is isomorphic to some $g^rD^{s,t}$ with $s+t\leq n-1$ by
corollary \ref{isom:cor}.
\item Note that in general ${\rm Ker}\, c$ is not irreducible, see
\cite{Kulk} for the  reduction matter.
\end{enumerate}
\end{remark}
\subsection{Orthogonal decomposition of $D^{p,q}$}
Following Kulkarni we call the elements in $\ker c\subset D^{p,q}$
effective elements of $D^{p,q}$. And shall be denoted by $E^{p,q}$.
\par\noindent
So if we apply corollary \ref{orth:decomp} several times we obtain the orthogonal decomposition
of $D^{p,q}$:
\begin{equation}\label{ort:decom}
D^{p,q}=E^{p,q}\oplus gE^{p-1,q-1}\oplus g^2E^{p-2,q-2}\oplus ...
\oplus g^rD^{p-r,q-r}
\end{equation}
where $r={\min \{p,q\}}$.\par\noindent
In this section, we will see how double forms decompose explicitly under this
orthogonal decomposition. To simplify the exposition, we shall consider only
the case where $p=q$.
\par
First, note that formula \ref{ck:gl}, for $\omega \in E^{p,p}$ becomes
\begin{equation}
\begin{split}
\label{ck:glef}
c^k({g^l\over l!}\omega)=\prod_{i=1}^{i=k}(n-2p-l+i){g^{l-k}\over (l-k)!}.
\omega\qquad &{\rm if}\quad l\geq k\\
c^k(g^l\omega)=0\qquad &{\rm if}\quad l<k
\end{split}
\end{equation}
With respect to the previous orthogonal  decomposition,
let $\omega=\sum_{i=0}^{p}g^i\omega_{p-i}\in D^{p,p}$ where
$\omega_{p-i}\in E^{p-i,p-i}$, then using the previous formula
\ref{ck:glef},
 we have
\begin{equation*}
\begin{split}
c^k(\omega)&=\sum_{i=0}^{p}c^k(g^i\omega_{p-i})=\sum_{i=k}^{p}c^k(g^i\omega_{p-i})\\
           &=\sum_{i=k}^{p}i!\prod_{j=1}^{j=k}(n-2(p-i)-i+j)
           {g^{i-k}\over (i-k)!}\omega_{p-i}
\end{split}
\end{equation*}
Therefore, we get
\begin{equation}\label{kom:bom}
c^k(\omega)
= \sum_{i=k}^{p}i! \prod_{j=1}^{j=k}(n-2p+i+j) {g^{i-k}\over (i-k)!}\omega_{p-i}
\end{equation}
Taking in the previous formula $k=p,p-1,p-2,..,k,..,0$ respectively,
 and solving for $\omega_k$ we get
\begin{equation*}
\begin{split}
{p!n!\over (n-p)!}\omega_0&= c^p(\omega)\\
{(p-1)!(n-2)!\over (n-p-1)!}\omega_1&= c^{p-1}(\omega)-
{1\over n}g.c^p(\omega)\\
{(p-2)!(n-4)!\over (n-p-2)!}\omega_2&= c^{p-2}(\omega)-
{1\over n-2}gc^{p-1}(\omega)+
{1\over 2!(n-2)(n-1)}g^2c^{p}(\omega)\\
&\vdots\\
{(p-k)!(n-2k)!\over (n-p-k)!}\omega_k&=c^{p-k}(\omega)+
\sum_{r=1}^k{(-1)^r\over {r!\prod_{i=0}^{r-1}(n-2k+2+i)}}
g^rc^{p-k+r}(\omega)\\
&\vdots\\
\omega_p&=\omega+\sum_{r=1}^p{(-1)^r\over {r!\prod_{i=0}^{r-1}
(n-2p+2+i)}}g^rc^{r}(\omega)
\end{split}
\end{equation*}
Note that $\omega_p={\rm con}\, \omega$ is the conformal component
 defined by Kulkarni.
\par\noindent
We have therefore proved the following theorem which generalizes a similar classical result in the case where $\omega$ is the Riemann curvature tensor:\par\smallskip\noindent
\begin{theorem}\label{omk:form} With respect to the orthogonal decomposition \ref{ort:decom},
 each  $\omega\in D^{p,p}$ is decomposed as follows
$$\omega=\omega_p+g.\omega_{p-1}+...+g^p.\omega_0$$
where
\begin{equation*}
\omega_k={(n-p-k)!\over (p-k)!(n-2k)!}\left[ c^{p-k}(\omega)+
\sum_{r=1}^k{(-1)^r\over {\prod_{i=0}^{r-1}(n-2k+2+i)}}
{g^r\over r!}c^{p-k+r}(\omega)\right]
\end{equation*}
\end{theorem}
In particular, for $\omega=R$, we recover the well known decomposition of Riemann curvature tensor
$$R=W+{1\over n-2}(c(R)-{1\over n}g.c^2(R))g+{1\over 2n(n-1)}c^2(R).g^2$$
\section{The algebra of curvature structures}
Remark that from the definition of the product (see formula \ref{def:prod}), we have
 $$\omega_1.\omega_2=(-1)^{pr+qs}\omega_2.\omega_1$$
 Then following,
 Kulkarni, we define the algebra of curvature structures
 to be the commutative sub-algebra ${\cal C}=\bigoplus_{p\geq 0}{\cal C}^p$,
  where ${\cal C}^p$ denotes the symmetric
 elements of ${D}^{p,p}$. That is the sub-algebra of
  symmetric double forms.\par\medskip
  Another basic map in ${D}^{p,q}$ is the first Bianchi sum,
  denoted ${\cal B}$. It maps  ${\cal D}^{p,q}$ into
 ${\cal D}^{p+1,q-1}$ and is defined as follows. Let
 $\omega \in  {\cal D}^{p,q}$, set ${\cal B}\omega=0$ if $q=0$.
 Otherwise  set
$${\cal B}\omega(x_1\wedge...\wedge x_{p+1},y_1\wedge...\wedge y_{q-1})=\sum_{j=1}^{p+1}(-1)^j\omega(x_1\wedge...\wedge \hat{x}_j\wedge ... x_{p+1},x_j\wedge y_1\wedge...\wedge y_{q-1})$$
where $\hat{}$ denotes omission.\par\smallskip\noindent
It is easy to show that for $\omega\in{\cal D}^{p,q},\theta\in{\cal D}^{r,s}$, one have \cite{Kulk}
$${\cal B}(\omega.\theta)={\cal B}\omega.\theta+(-1)^{p+q}\omega.{\cal B}\theta$$
Consequently, ${\rm ker}{\cal B}$ is closed under multiplication in  ${\cal D}$.\par\smallskip\noindent
 The algebra of curvature structures satisfying the first Bianchi identity is
 defined to be  ${\cal C}_1={\cal C}\cap {\rm ker}{\cal B}$
\par\medskip\noindent
\subsection{Sectional curvature}
 Let $G_p$ denote the Grassman algebra of $p$-planes in $V$, and $\omega\in
 {\cal C}^p $. We define the sectional curvature of $\omega$   to be
 $$K_\omega(P)=\omega(e_1\wedge...\wedge e_p,e_1\wedge...\wedge e_p)$$
 where $\{e_1,...,e_p\}$ is any orthonormal basis of $V$. \par
 \smallskip\noindent
 Using formula \ref{eps:prod}, we can evaluate the sectional curvature of the tensors
 $g^p\omega$
 for $\omega \in {\cal C}^r$ and $\{e_1,...,e_{p+r}\}$ orthonormal , as follows
 \begin{equation}\label{gpom:form}
 \begin{split}
  g^p\omega&(e_1\wedge...\wedge e_{p+r},e_1\wedge...\wedge e_{p+r})\\
 =&p!\sum_{1\leq i_1<i_2<...<i_r\leq p+r}\omega(e_{i_1}\wedge...\wedge e_{i_r},
 e_{i_1}\wedge...\wedge e_{i_r})=p!{\rm trace}\,\omega_{\mid \Lambda^rP}
 \end{split}
 \end{equation}
 where $P$ denotes the plane spanned by $\{e_1,...,e_{p+r}\}$.
 \par\medskip\noindent
 The sectional curvature $K_{\omega}$ determines generically $\omega$.
 Precisely,
 for   $\omega, \theta
 \in {\cal C}_1^p$, the equality $K_{\omega}=K_{\theta}$ implies $\omega=\theta$
   (cf. prop. 2.1 in \cite{Kulk}). In particular we have the following characterization
   of the curvature structures $\omega \in {\cal C}_1^p$ with constant sectional curvature
  \begin{equation}\label{kom:const}
  K_{\omega}\equiv c \qquad \text{ if and only if}\qquad
  \omega=c{g^p\over p!}
  \end{equation}
Next, we shall prove a useful explicit formula for the Hodge star operator
\begin{theorem}
For $\omega\in {\cal C}_1^p$ and  $1\leq p\leq k\leq n$ we have
\begin{equation}\label{star:form}
{1\over (k-p)!}*(g^{k-p}\omega)=\sum_{r={\max\{0,p-n+k\}}}^p
{(-1)^{r+p}\over r!}{g^{n-k-p+r}\over (n-k-p+r)!}c^r\omega
\end{equation}
In particular for $k=n$ and $k=n-1$ respectively, we get
\begin{equation}\label{nn-1:star}
*({g^{n-p}\omega\over (n-p)!})={1\over p!}c^p\omega\quad {\rm and}\,
*({g^{n-p-1}\omega\over (n-p-1)!})={c^p\omega\over p!} g-
{c^{p-1}\over (p-1)!}\omega
\end{equation}
\end{theorem}
{\bf Proof.} It is not difficult to check that
$$\sum_{1\leq i_1,i_2,...,i_p\leq n}\omega(e_{i_1}\wedge ...\wedge e_{i_p},e_{i_1}\wedge ...\wedge e_{i_p})= \sum_{1\leq i_1,i_2,...,i_p\leq k}\omega(e_{i_1}\wedge ...\wedge e_{i_p},e_{i_1}\wedge ...\wedge e_{i_p})$$
$$+\sum_{r=0}^{p-1}(-1)^{r+p+1}C_r^p \sum_{k+1\leq i_{r+1},...,i_p\leq n}c^r\omega(e_{i_{r+1}}\wedge ...\wedge e_{i_p},e_{i_{r+1}}\wedge ...\wedge e_{i_p})$$
Then using formula \ref{gpom:form} the previous formula becomes
\begin{equation*}
\begin{split}
c^p\omega=p! &{g^{k-p}\over (k-p)!}\omega(e_{i_1}\wedge ...\wedge e_{i_k},
e_{i_1}\wedge ...\wedge e_{i_k})+
\sum_{r={\max\{0,p-n+k\}}}^{p-1} (-1)^{r+p+1}  C_r^p \\
&{(p-r)! \over (n-k-p+r)!}g^{n-k-p+r}c^r\omega
(e_{i_{k+1}}\wedge ...\wedge e_{i_n},e_{i_{k+1}}\wedge ...\wedge e_{i_n})
\end{split}
\end{equation*}
Finally, note that the general term of the previous sum is $c^p\omega$ if
 $r=p$. This completes the proof, since both sides of the equation satisfy the first Bianchi
 identity.\ppp
 \par\smallskip\noindent
 The following corollary is a direct consequence of the previous theorem.
\begin{corollary}\label{starformula}\begin{enumerate}
 \item
For $\omega\in {\cal C}_1^p$ and  $1\leq p\leq n$ we have
\begin{equation}
*\omega=\sum_{r=\max\{0,2p-n\}}^p
{(-1)^{r+p}\over r!}{g^{n-2p+r}\over (n-2p+r)!}c^r\omega
\end{equation}
\item
For all $0\leq k\leq n$ we have
$$*{g^k\over k!}={g^{n-k}\over (n-k)!}$$
\end{enumerate}
\end{corollary}
\begin{theorem} With respect to the decomposition \ref{omk:form},
 we have for $\omega=\sum_{i=0}^p g^{p-i}\omega_i$
\begin{equation}\label{star:omi}
*\omega=\sum_{i=0}^{\min\{p,n-p\}}(p-i)!(-1)^i
{1\over (n-p-i)!}g^{n-p-i}\omega_i
\end{equation}
In particular if $n=2p$, we have
$$*\omega=\sum_{i=0}^p(-1)^ig^{p-i}\omega_i$$
\end{theorem}
{\bf Proof.} First, let $\omega\in E_1^{i}$ be effective then formula
 \ref{star:form} shows that
\begin{equation}
{1\over (k-i)!}*(g^{k-i}\omega)=
\begin{cases}
0&\text{if $i-n+k>0$},\\
{(-1)^i\over (n-k-i)!}g^{n-k-i}\omega &\text{if $i-n+k\leq 0$}
\end{cases}
\end{equation}
Next, let $\omega=\sum_{i=0}^p g^{p-i}\omega_i$, where $\omega_i\in E_1^i$,
then
\begin{equation*}
\begin{split}
*\omega=& \sum_{i=0}^p *(g^{p-i}\omega_i)\\
=& \sum_{i=0}^{\min\{p,n-p\}} (p-i)! {(-1)^i\over (n-p-i)!}g^{n-p-i}\omega_i
\end{split}
\end{equation*}
\ppp
\begin{corollary} With respect to the decomposition \ref{ort:decom},
 we have for $\omega=\sum_{i=0}^p g^{p-i}\omega_i$
\begin{equation}
\label{stargl:omi}
*(g^l\omega)=\sum_{i=0}^{\min\{p,n-p-l\}} (p-i+l)! {(-1)^i\over (n-p-l-i)!}
g^{n-p-l-i}\omega_i
\end{equation}
\end{corollary}
{\bf Proof.} First  formula \ref{star:omi} implies that
$$*(g^l\omega)=\sum_{i=0}^{\min\{p+l,n-p-l\}} (p-i+l)! {(-1)^i\over (n-p-l-i)!}
g^{n-p-l-i}(g^l\omega)_i $$
Next, note that
$$g^l\omega=\sum_{i=0}^{p}g^{p+l-i}\omega_i$$
Consequently
\begin{equation*}
(g^l\omega)_i=
\begin{cases}
0 &\text {if $i>p$}\\
\omega_i&\text{if $i\leq p$}
\end{cases}
\end{equation*}
This completes the proof of the corollary.\ppp
\section{The $(p,q)$-curvatures}
Let $(M,g)$ be an $n$-dimensional Riemannian manifold and $T_mM$ be its tangent
space at a point $m\in M$. Let $D^{p,q},{\cal C}^p,{\cal C}_1^p...$
denote also the vector bundles over $M$ having as fibers at $m$, the spaces
$D^{p,q}(T_mM)$, ${\cal C}^p(T_mM)$, ${\cal C}_1^p(T_mM) ...$.
 Note that all the above  algebraic results, can be applied to the ring of
all global sections of these bundles.\par\medskip\noindent
Remark that since the metric $g$ and the Riemann curvature tensor $R$ both satisfy  the first Bianchi identity then so are
all the tensors $g^pR^q$ and $*(g^pR^q)$. The aim of this section is to sudy some geometric properties of these tensors.
First we start with the case $q=1$:
\subsection{The $p$-curvature}
Recall that the $p$-curvature \cite{Lab2,Lab3}, defined for $0\leq p\leq n-2$
 and denoted by
$s_p$,
 is the sectional curvature of the tensor
$${1\over (n-2-p)!}*(g^{n-2-p}R)$$
For a given tangent $p$-plane at $m\in M$, $s_p(P)$ coincides with the half of
the scalar curvature at $m$
of the totally geodesic submanifold ${\rm exp}_m\, P^{\bot}$. For $p=0$ it is
the the half of the usual scalar curvature, and for $p=n-2$ it coincides with the usual sectional
curvature.\par\medskip\noindent
 In this subsection, using the $p$-curvature and the previous results,
we shall give a short proof for the following properties. Similar results
were proved by a long calculation in \cite{Shin} and
\cite{Lab2}.
\begin{theorem}\label{sp+sp:sp-sp}\begin{enumerate}
\item For each  $2\leq p\leq n-2$, the $p$-curvature is
constant if and only if $(M,g)$ is with constant sectional curvature.
\item For each $1\leq p\leq n-1$,  the Riemannian manifold $(M,g)$ is
 Einstein if and only if the function $P\rightarrow s_p(P)-s_{n-p}(P^{\bot})
 =\lambda$ is
 constant. Furthermore, in such case we have $\lambda={n-2p\over 2n}c^2R.$
\item For each  $2\leq p\leq n-2$ and $p\not = {n\over 2}$, the function
  $P\rightarrow s_p(P)+s_{n-p}(P^{\bot})=\lambda$ is constant if and only if
 the  manifold $(M,g)$ is with constant sectional curvature. Furthermore,
  in such case we have $\lambda={2p(p-1)+(n-2p)(n-1)\over 2n(n-1)}c^2R$
\item Let $n=2p$. Then  the Riemannian
 manifold $(M,g)$ is  conformally flat with constant scalar curvature
  if and only if the function
  $P\rightarrow s_p(P)+s_{p}(P^{\bot})=\lambda$ is constant. Furthermore,
  in such case we have $\lambda={n-2\over4(n-1)} c^2R.$
\end{enumerate}
\end{theorem}
 {\bf Proof.} First we prove 1).  Let $s_p \equiv c$ then the sectional
curvature of the tensor $g^{n-2-p}R\in {\cal C}_1^{n-p}$ is constant.
 Therefore we have ${1\over (n-2-p)!}g^{n-2-p}R=c{g^{n-p}\over (n-p)!}$
  and so by proposition \ref{mult:injec} we have $(n-p)(n-p-1)R=cg^2$. That is $R$ is with
   constant sectional curvature. \par
Next we prove 2). Suppose $s_p(P)-s_{n-p}(P^{\bot})=c$ for all $P$, then
$${1\over (n-2-p)!}g^{n-2-p}R(*P,*P)-{1\over (p-2)!}g^{p-2}R(P,P)=c \, \text{for all P}$$
then using formula \ref{star:form} we get
$$\sum_{r=0}^2{(-1)^{r}\over r!}{g^{p-2+r}\over (p-2+r)!}c^rR(P,P)-{1\over (p-2)!}g^{p-2}R(P,P)=c \, \text{for all P}$$
The left hand side is the sectional curvature of a curvature tensor which satisfies the first Bianchi identity then
$$-{g^{p-1}\over (p-1)!}cR+{g^{p}\over 2(p!)}c^2R=c{1\over p!}g^{p}$$
proposition \ref{mult:injec}, implies that
$$-cR+{g\over 2p}c^2R=c{1\over p}g$$
and therefore
$$cR={c^2R-2c\over 2p}g$$
so that $(M,g)$ is an Einstein manifold. Furthermore, after taking the trace we get $c={n-2p\over
2n}c^2R$.\par
 Finally we prove 3) and 4). Suppose $s_p(P)+s_{n-p}(P^{\bot})=c$ for all $P$,
  then as in part  2) we have,
$$2{g^{p-2}\over (p-2)!}R-{g^{p-1}\over (p-1)!}cR+{g^{p}\over 2(p!)}c^2R=c{1\over p!}g^{p}$$
then using \ref{mult:injec} we get
$$2R-{g\over (p-1)}cR+{g^{2}\over 2p(p-1)}c^2R=c{1\over p(p-1)}g^{2}$$
Which  implies that
$$2\omega_2-{n-2p\over p-1}g\omega_1+(2+{(n-2p)(n-1)\over p(p-1)})g^2\omega_0
={c\over p(p-1)}g^2$$
where $R=\omega_2+g\omega_1+g^2\omega_0$. Then if $n\not =2p$
then $\omega_2=\omega_1=0$ and therefore the sectional curvature of $(M,g)$
 is constant. In the case $n=2p$, we have $\omega_2=0$ and
 $c=2\omega_0p(p-1)={n-2\over 4(n-1)}c^2R$. So that $(M,g)$
  is conformally flat with constant scalar curvature.  \ppp

 \subsection{The $(p,q)$-curvatures}
 The $(p,q)$-curvatures are the $p$-curvatures of the Gauss-Kronecker tensor
 $R^q$ (that is the product of the Riemann tensor R with itself k-times
  in the ring of
 curvature structures). Precisely, they are defined by\par\noindent
\begin{definition} The $(p,q)$-curvature, denoted $s_{(p,q)}$,
  for $1\leq q\leq {n\over 2}$ and $0\leq p\leq n-2q$, is the sectional curvature
  of the following $(p,q)$-curvature tensor
  \begin{equation}\label{spq:def}
  R_{(p,q)}={1\over (n-2q-p)!}*\bigl( g^{n-2q-p}R^q\bigr)
  \end{equation}
 In other words, $s_{(p,q)}(P)$ is the sectional curvature of the tensor
  ${1\over (n-2q-p)!}g^{n-2q-p}R^q$ at the orthogonal complement
of $P$.
\end{definition}
These curvatures include many of the well known curvatures.\\
Note that for $q=1$, we have  $s_{(p,1)}=s_p$ coincides with
 the $p$-curvature. In particular,  $s_{(0,1)}$ is the half of
  the scalar curvature and  $s_{(n-2,1)}$ is
   the sectional curvature of $(M,g)$.\\
For $p=0$ and $2q=n$, $s_{(0,{n\over 2})}=*R^{n/2}$ is up to a constant
 the Killing-Lipshitz curvature. More generally, $s_{(n-2q,q)}(P)$ is,
  up to a constant, the Killing-Lipshitz curvature
of $P^{\bot}$. That is the $(2p)$-sectional curvatures defined by A. Thorpe
in \cite{Thorpe}.\par
For $p=0$, $s_{(0,q)}=*{1\over (n-2q)!}g^{n-2q}R^q={1\over (2q)!}c^{2q}R^q$
 are scalar functions which generalize the usual scalar curvature.
 They are up to  constants the
 integrands in the Weyl tube formula \cite{Wey}.
\par\noindent
For $p=1$, $s_{(1,q)}$ are the curvatures of
 generalized Einstein tensors. Precisely, let us define the following:\par\smallskip\noindent
\begin{definition}\begin{enumerate}
\item The $2q$-scalar curvature function, or the $2q$-H. Weyl curvature invariant,  denoted $h_{2q}$, is
the $(0,q)$-curvature. That is  $$h_{2q}=s_{(0,q)}={1\over (2q)!}c^{2q}R^q$$
\item The $2q$-Einstein tensor, denoted $T_{2q}$, is defined to be the $(1,q)$-curvature
tensor, that is
$$T_{2q}=*{1\over (n-2q-1)!}g^{n-2q-1}R^q$$
\end{enumerate}
\end{definition}
By formula \ref{nn-1:star} we have
$$T_{2q}={1\over (2q)!}c^{2q}R^q-{1\over (2q-1)!}c^{2q-1}R^q=h_{2q}-
{1\over (2q-1)!}c^{2q-1}R^q$$
For $q=1$ we recover the usual Einstein tensor $T_2={1\over 2}c^2R-cR$.
Note that $c^{2q-1}R^q$ can be considered as a generalization of the Ricci curvature.
In a forthcoming paper \cite{Lab3} we shall prove that the $2q$-Einstein tensor is the
gradiant of the total $2q$-scalar curvature function seen as a functional on the
space of all Riemannian metrics with volume 1.
Which generalize the classical well known result about the scalar curvature.\\
Finally note that in general $s_{(p,q)}(P)$ coincides also with the
$2q$-scalar curvature of $P^{\bot}$. \par\medskip\noindent
\subsection{Examples}
1) Let $(M,g)$ be with constant sectional curvature $\lambda$, then
$$R={\lambda\over 2}g^2\qquad {\text and}\qquad R^q={\lambda^q\over 2^q}g^{2q}$$
And therefore
$$*{1\over (n-2q-p)!}g^{n-2q-p}R^q=*{\lambda^q\over 2^q(n-2q-p)!}g^{n-p}
={\lambda^q(n-p)!\over 2^q(n-2q-p)!}{g^p\over p!}$$
so that  the $(p,q)$-curvature is also  constant  and equal to
${\lambda^q(n-p)!\over 2^q(n-2q-p)!}$. \par\noindent
The converse will be discussed in the next section.\par\bigskip\noindent
2) Let $(M,g)$ be a Riemannian product of two Riemannian manifolds
 $(M_1,g_1)$ and $(M_2,g_2)$. If we index by $i$ the invariants of the
 metric $g_i$ for $i=1,2$, then
 $$R=R_1+R_2\quad {\rm and}\quad R^q=(R_1+R_2)^q=
 \sum_{i=0}^qC_i^qR_1^iR_2^{q-i}$$
 consequently, a straightforward calculation shows that
 \begin{equation*}
 \begin{split}
 h_{2q}&={c^{2q}R^q\over (2q)!}=\sum_{i=0}^qC_i^q {c^{2q}\over (2q)!}
 (R_1^iR_2^{q-i})\\
 &=\sum_{i=0}^qC_i^q {c^{2i}R_1^i\over (2i)!}
 {c^{2q-2i}R_2^{q-i}\over (2q-2i)!}\\
 &=\sum_{i=0}^qC_i^q (h_{2i})_1 (h_{2q-2i})_2
\end{split} \end{equation*}
Where we used the convention $h_0=1$.
\par\bigskip\noindent
3) Let $(M,g)$  be a hypersurface of the Euclidean space.  If $B$ denotes the second
fundamental form at a given point, then the Gauss equation shows that
$$R={1\over 2}B^2\qquad {\text and}\quad R^q={1\over 2^q}B^{2q}$$
Consequently, if $\lambda_1\leq \lambda_2\leq ...\leq \lambda_n$ denote the eigenvalues
of $B$, then the eigenvalues of $R^q$ are ${(2q)!\over 2^q}\lambda_{i_1}
\lambda_{i_2}...\lambda_{i_{2q}}$ where $ i_1<...<i_{2q}$.
Therefore all the tensors $g^pR^q$ are diagonalizable and their eigenvalues
have the following form
$$g^pR^q(e_1... e_{p+2q},e_1... e_{p+2q})=
 {p!(2q)!\over 2^q}\sum_{1\leq i_1<...<i_{2q} \leq p+2q}\lambda_{i_1}
 ...\lambda_{i_{2q}}$$
 where $\{e_1,... e_n\}$ is an orthonormal basis of eigenvectors of $B$.
 In particular we have
$$h_{2q}=s_{(0,q)}=
 {(2q)!\over 2^q}\sum_{1\leq i_1<...<i_{2q} \leq n}\lambda_{i_1}
 ...\lambda_{i_{2q}}$$
 So they are, up to a constant, the symmetric functions in the eigenvalues
 of $B$.\par\noindent
 and more generally, we have
$$s_{(p,q)}(e_{n-p+1},...,e_n)=
 {(2q)!\over 2^q}\sum_{1\leq i_1<...<i_{2q} \leq n-p}\lambda_{i_1}
 ...\lambda_{i_{2q}}$$
 \par\medskip\noindent
 4)  Let $(M,g)$  be a  conformally flat manifold. Then it is well known that
 at each point
 of $M$,
 the Riemann curvature tensor is determined by a symmetric bilinear
  form $h$,
 precisely we have $R=g.h$. Consequently, $R^q=g^qh^q$.\par
 Let $\{e_1,... e_n\}$ be an orthonormal basis of eigenvectors of $h$ and
 $\lambda_1\leq \lambda_2\leq ...\leq \lambda_n$ denote the eigenvalues of $h$.
 \par
 Then it is not difficult to see that all the tensors $g^pR^q$ are also
  diagonalizable.
 The eigenvalues are given by
$$g^pR^q(e_1... e_{p+2q},e_1... e_{p+2q})=
 (p+q)!q!\sum_{1\leq i_1<...<i_{q} \leq p+2q}\lambda_{i_1}
 ...\lambda_{i_{q}}$$
 In particular the $(p,q)$-curvatures are determined by
\begin{equation*}
s_{(p,q)}(e_{n-p+1},...,e_n)={(n-q-p)!q!\over (n-2q-p)!}
\sum_{1\leq i_1<...<i_{q} \leq n-p}\lambda_{i_1}
 ...\lambda_{i_{q}}
 \end{equation*}
\subsection{Properties}
 The following theorem generalizes a similar induction formula \cite{Lab2}
  for
 the $p$-curvature:\par\medskip\noindent
 \begin{theorem} For $1\leq q\leq {n\over 2}$ and
  $1\leq p\leq n-2q$   we have
  $$\sum_{k=p}^n s_{(p,q)}(P,e_k)=(n-2q-p+1)s_{(p-1,q)}(P)$$
  where $P$ is an arbitrary tangent $(p-1)$-plane and $\{e_{p},...,e_n\}$
  is any orthonormal basis of $P^\bot$.\\
  In particular we have
  $$\sum_{i=1}^n T_{2q}(e_i,e_i)=(n-2q)h_{2q}$$
  \end{theorem}
  {\bf Proof.} Using \ref{gstar:c} we have
  \begin{equation*}
  \begin{split}
  {1 \over (n-2q-p)!}c*(g^{n-2q-p}R^q)&={1 \over (n-2q-p)!}*g(g^{n-2q-p}R^q)\\
&= (n-2q-p+1)*({g^{n-2q-p+1}\over (n-2q-p+1)!} R^q)
\end{split}
\end{equation*}
  to finish the proof just take the sectionnal curvatures of both
  sides.\ppp
 \par\medskip\noindent
 The following proposition is the only exception in this paper where
 one
 needs
 the use of the second Bianchi identity, see \cite{Lab3} for the proof:
 \begin{proposition}[Schur's theorem]
 Let $p\geq 1$ and $q\geq 1$.
 If at every point $m\in M$ the $(p,q)$-curvature is constant
 (that is on the fiber at m),
 then it is constant.\end{proposition}
 The following can be seen as the converse of a Thorpe's result \cite{Thorpe}
 \begin{proposition}
 If $R^s$ and $R^{s+r}$ are both with constant sectional curvature $\lambda$
 and $\mu$ respectively, such that $\lambda \not=0$ and $s+2r\leq n$. Then
  $R^r$ is also  with constant sectional curvature and equal to
  ${\mu s!r!\over\lambda(s+r)!}$.\end{proposition}
 {\bf Proof.} Suppose
 $$R^s=\lambda{g^s\over s!}\qquad {\rm and}
\quad R^{s+r}=\mu{g^{(s+r)}\over (s+r)!}$$
 then
 $$\lambda{g^s\over s!}R^r=\mu{g^{(s+r)}\over (s+r)!}$$
 then since $s+2r\leq n$, proposition \ref{mult:injec} shows that
 $${\lambda\over s!}R^r=\mu{g^r\over (s+r)!}$$
 This completes the proof.\ppp
 \par\noindent\medskip
 Like the case of $R^s$, it is not true in general that if $h_{2s}$ is
 constant then the higher
 scalar curvatures are constant, nevertheless we have the following result
 \begin{proposition} If for some $s$ the tensor $R^s$ is with constant
 sectionnal curvature $\lambda$ then we have for all $r\geq o$
 $$h_{2s+2r}={(n-2r)!\over (2s)!(n-2s-2r)!}\lambda h_{2r}$$
 In particular, if $n$ is even the Gauss-Bonnet integrand is determined by
 $$h_n=\lambda h_{n-2s}$$
 \end{proposition}
 {\bf Proof.} Suppose $R^s=\lambda{g^{2s}\over (2s)!}$ then
 \begin{equation*}
 \begin{split}
 h_{2s+2r}=& {1\over (n-2s-2r)!}*(g^{n-2s-2r}R^{s+r})\\
 =& {1\over (n-2s-2r)!}*(g^{n-2s-2r}\lambda{g^{2s}\over (2s)!}R^r )\\
 =&  {\lambda\over (2s)! (n-2s-2r)!}*(g^{n-2r}R^r )=
 {(n-2r)!\lambda\over (2s)! (n-2s-2r)!}h_{2r}
 \end{split}
 \end{equation*}
 \ppp
 \begin{theorem} \begin{enumerate}
 \item For every $(p,q)$ such that $2q\leq p\leq n-2q$,
 the $(p,q)$-curvature
 $s_{(p,q)}\equiv \lambda$ is constant if and only if the sectional curvature of
 $R^q$ is constant and equal to ${\lambda(2q)!(n-p-2q)!\over (n-p)!}$.
 \item For every $(p,q)$ such that $p<2q$, the $(p,q)$-curvature
 $s_{(p,q)}\equiv c$ is constant if and only if $c^{2q-p}(R^q)$ is proportional
 to the metric. That is $c^{2q-p}(R^q)=const.g^p$.
 \end{enumerate}
 \end{theorem}
{\bf Proof.} Recall that  $s_{(p,q)}\equiv \lambda$ if and only if
$$ {g^{n-2q-p}\over (n-2q-p)!}R^q=\lambda {g^{n-p}\over (n-p)!}$$
that is
$$ g^{n-2q-p}\bigl({R^q\over (n-2q-p)!}-\lambda {g^{2q}\over (n-p)!}\bigr)=0$$
 now, Let $2q\leq p\leq n-2q$, then by proposition \ref{mult:injec}
 this is equivalent to
$$R^q=\lambda{(n-2q-p)!\over (n-p)!}g^{2q}$$
Next, if $p<2q$, then by remark 2.1 in section 2, our condition
 is equivalent
to
$$ c^{2q-p}\bigl({R^q\over (n-2q-p)!}-\lambda {g^{2q}\over (n-p)!}\bigr)=0$$
that is
$$ c^{2q-p}(R^q)=const. g^p$$
which completes the proof of the theorem.\ppp
\par\medskip\noindent
The following lemma provides a characterization of the previous condition
on $R^q$ and generalizes a similar result in the case
 of Ricci curvature ($p=q=1$).
\begin{lemma} For $p<2q$, the tensor $c^{2q-p}(R^q)$ is proportional to the metric $g^p$
if and only if
$$\omega_i=0 \, {\text { for}}\,  1\leq i\leq {\min\{p,n-p\}}$$
where $R^q=\sum_{i=0}^{2q}g^{2q-i}\omega_i$.
\end{lemma}
{\bf Proof.} Formula \ref{kom:bom} shows that
$$ c^{2q-p}(R^q)=\sum_{i=2q-p}^{2q}i!
\biggl(\prod_{j=1}^{2q-p}(n-4q+i+j)\biggr)
{g^{i-2q+p}\over (i-2q+p)!}\omega_{2q-i}$$
and therefore $c^{2q-p}(R^q)=\lambda g^p$ if and only if
$$\sum_{s=0}^p(2q-s)! \biggl(\prod_{j=1}^{2q-p}(n-2q-s+j)\biggr)
{g^{p-s}\over (p-s)!}\omega_{s}=\lambda g^p$$
where we changed the index to $s=2q-i$. consequently,
$$g^{p-s}\omega_s=0\, {\text {for}}\,  1\leq s\leq p, {\text {and}}\,
 \lambda={(2q)!\over p!}\bigl(\prod_{j=1}^{2q-p}(n-2q+j)\bigr)\omega_0$$
 By proposition \ref{mult:injec}, this is equivalent to $\omega_s=0$ for $1\leq s\leq n-p$ and $s\leq p$.
 Note that $ g^{p-s}\omega_s=0$ if $s>n-p$.
 This completes the proof of lemma.\ppp
\begin {theorem}\label{pqcurv:thm}\begin{enumerate}
\item Let $2q\leq r\leq n-2q$, $n\not =2r$ and
$R^q=\sum_{i=0}^{2q}g^{2q-i}\omega_i$,
then\begin{enumerate}
\item The function $P\rightarrow s_{(r,q)}(P)- s_{(n-r,q)}(P^\bot) \equiv
 \lambda$ is constant
if and only if $\omega_i=0$ for $1\leq i\leq 2q-1$ and
 $\bigl({(n-r)!\over (n-2q-r)!}-{r!\over (r-2q)!}\bigr)\omega_0=\lambda$.
\item The function $P\rightarrow s_{(r,q)}(P)+ s_{(n-r,q)}(P^\bot) \equiv
 \lambda$ is constant
if and only if $\omega_i=0$ for $1\leq i\leq 2q$ and $\bigl({(n-r)!\over
 (n-2q-r)!}+{r!\over (r-2q)!}\bigr)\omega_0=\lambda$. That is $R^q$ is with
  constant sectional curvature.\end{enumerate}
\item  Let $2q\leq r\leq n-2q$ and $n =2r$, then \begin{enumerate}
\item The function $P\rightarrow s_{(r,q)}(P)- s_{(r,q)}(P^\bot) \equiv
\lambda$ is constant
if and only if $\omega_i=0$ for $i$ odd such that $1\leq i\leq 2q-1$
 and $\lambda=0$
\item The function $P\rightarrow s_{(r,q)}(P)+ s_{(r,q)}(P^\bot) \equiv
 \lambda$ is constant
if and only if $\omega_i=0$ for $i$ even and $2\leq i\leq 2q$ and
 $2{r!\over (r-2q)!}\omega_0=\lambda$.
 \end{enumerate}
 \end{enumerate}
 \end{theorem}
{\bf Proof.}
Let $k,l\geq 0$ be such that $k+p=n-l-p$ and $\omega=\sum_{i=0}^p g^{p-i}\omega_i\in {\cal C}^p$, then

$$
{g^k\over k!} \omega -*({g^l\over l!}\omega)=
{g^k\over k!}\omega -\sum_{i=0}^{\min\{p,n-p-l\}}
{(p-i+l)!(-1)^i\over l!(n-p-l-i)!}g^{n-p-l-i}\omega_i$$
$$={1\over (n-2p-l)!} \sum_{i=0}^p \bigl[1-(-1)^i{(p-i+l)! (n-2p-l)!\over
 l!(n-p-l-i)!}\bigr]g^{n-p-l-i}\omega_i$$
therefore,
\begin{equation}\label {eq:1}
{g^k\over k!} \omega -*({g^l\over l!}\omega)=
\lambda{g^{n-l-p}\over (n-l-p)!}
\end{equation}
 if and only if
$$ \sum_{i=1}^p \bigl[1-(-1)^i{(p-i+l)! (n-2p-l)!\over  l!(n-p-l-i)!}\bigr]
g^{n-p-l-i}\omega_i\phantom{mmmmmm}$$
$${\phantom{ mmmm}}
+\bigl[\omega_0-{(p+l)! (n-2p-l)!\over  l!(n-p-l)!}\omega_0
- \lambda{(n-2p-l)!\over (n-l-p)!}\bigr]g^{n-l-p}=0$$
For $1\leq i\leq p$, let
$${\alpha}_i=1-(-1)^i{(p-i+l)! (n-2p-l)!\over  l!(n-p-l-i)!}
=1-(-1)^i{(s+l)!\over l!}{k!\over (s+k)!}$$
 where $s=p-i\leq p-1$.
 It is clear that $\alpha_{2j+1}>0$ and
 it is  not difficult to check that $\alpha_i\not=0$ for $i$ even,
 $1\leq i\leq p-1$ and   $k\not = l$. Also note that $\alpha_p= 1-(-1)^p=0$
 since p is even.\par
 Therefore in the case where  $k\not = l$, the condition \ref{eq:1}
  is equivalent to
 $$\omega_i=0\, {\text for}\, 1\leq i\leq p-1 \text{and}\,
 \lambda=\{{(n-p-l)!\over (n-2p-l)!}-{(p+l)!\over l!}\}\omega_0$$
 In the case $k=l$, we have $\alpha_i=1-(-1)^i$ the condition \ref{eq:1}
  is therefore
 equivalent to
 $$\omega_i=0\quad {\text {for $i$ odd,} }\quad 1\leq i\leq p \quad{\text {and} }\quad
  \lambda=0$$
  In a similar way, we have
 \begin{equation}\label{eq:2}
 {g^k\over k!} \omega +*({g^l\over l!}\omega)=
\lambda{g^{n-l-p}\over (n-l-p)!}
\end{equation}
 if and only if
\begin{equation*}
\begin{split}
 \sum_{i=1}^p &\bigl[1+(-1)^i{(p-i+l)! (n-2p-l)!\over  l!(n-p-l-i)!}\bigr]
g^{n-p-l-i}\omega_i\\
&{\phantom{mmmm}}+\bigl[\omega_0+{(p+l)! (n-2p-l)!\over
 l!(n-p-l)!}\omega_0
- \lambda{(n-2p-l)!\over (n-l-p)!}\bigr]g^{n-l-p}=0
\end{split}
\end{equation*}
In the case $k\not= l$ this is equivalent to
 $$\omega_i=0\; {\text {for}}\; 1\leq i\leq p\; {\text{and}}\;
  \lambda=\{{(n-p-l)!\over(n-2p-l)!} +
 {(p+l)!\over l!}\}\omega_0$$
 And in case $k=l$, we have ${\alpha}_i=1+(-1)^i$ and then condition
 \ref{eq:2} is then
 equivalent to
 $$\omega_i=0\quad {\text {for $i$ even,} }\quad 1\leq i\leq p\quad {\text {and} }\quad
  \lambda=2{(p+l)!\over l!}\omega_0$$
To complete the proof of the theorem just notice that, for $2q\leq r\leq n-2q$,  the condition $s_{(r,q)}\pm s_{(n-r,q)}\equiv \lambda$ is equivalent to
$$*({g^{n-2q-r}\over (n-2q-r)!}R^q)\pm {g^{r-2q}\over (r-2q)!}R^q=\lambda {g^r\over r!}$$
and then apply the previous result after taking $l=n-2q-r$ and $k=r-2q$
 and $p=2q$. \ppp
\par\medskip\noindent
{\bf Remark.} If $r<2q$, then $s_{(n-r,q)}\equiv 0$. So that our condition
 becomes $s_{(r,q)}\equiv \lambda$ is constant, and such case
  was discussed above.
\section{Generalized Avez type formula}
The following theorem generalizes a  result due to Avez \cite{Avez} in
 the case
when $n=4$ and $\omega=\theta=R$.
\begin{theorem} Let $n=2p$ and $\omega, \theta\in {\cal C}_1^p$, then
$$*(\omega\theta)=\sum_{r=0}^p{(-1)^{r+p}\over (r!)^2}<c^r\omega,c^r\theta>$$
In particular if $n=4q$, then the Gauss-Bonnet integrand is
determined by
$$h_{4q}=\sum_{r=0}^{2q}{(-1)^{r}\over (r!)^2}|c^rR^q|^2$$
\end{theorem}
{\bf Proof.} Let  $\theta\in {\cal C}_1^p$ and $\omega\in {\cal C}_1^{n-p}$,
then using formula \ref{in:star} and  corollary \ref{starformula} we get
$$*(\omega\theta)=<\omega,*\theta>=\sum_{r={\max\{0,2p-n\}}}^p{(-1)^{r+p}\over
(r!)(n-2p+r)!}<\omega,g^{n-2p+r}c^r\theta>$$
To complete the proof just take $n=2p$ and use theorem \ref{theo:gc}.\ppp
\par\noindent\medskip
The following corollary is an alternative way to write the previous formula.
\begin{corollary} Let $n=2p$ and $\omega, \theta\in {\cal C}_1^p$, then
$$*(\omega\theta)=\sum_{r=0}^p{(-1)^{r+p}\over (r!)^2}<g^r\omega,g^r\theta>$$
\end{corollary}
{\bf Proof.} It is a direct consequence of the previous formula
and corollary \ref{ckgk:gkck}\ppp
\par\noindent\medskip
The following result is of the same type as the previous one
\begin{theorem} With respect to the decomposition \ref{ort:decom}, let
  $\omega=\sum_{i=0}^{n-p}g^{n-p-i}\omega_i\in {\cal C}_1^{n-p}$ and
 $\theta=\sum_{i=0}^p g^{p-i}\theta_i\in {\cal C}_1^p$ , then
$$*(\omega\theta)=\sum_{r=0}^{\min\{p,n-p\}}(-1)^r(n-2r)!<\omega_i,\theta_i>$$
\end{theorem}
{\bf Proof.} By formula \ref{star:omi} we have
$$*(\omega\theta)=<\omega,*\theta>=\sum^{{\min\{p,n-p\}}}_{i=0}{(-1)^i(p-i)!
\over(n-p-i)!}<\omega,g^{n-p-i}\theta_i>$$
and therefore using  lemma \ref{lemma:gpq} below, we get
$$ *(\omega\theta)=\sum^{{\min\{p,n-p\}}}_{i=0}{(-1)^i(p-i)!\over(n-p-i)!}
<g^{n-p-i}\omega_i,g^{n-p-i}\theta_i>$$
After considering separately the cases $p<n-p,p=n-p$ and $p>n-p$ and the lemma
below one can complete the proof easily.\ppp
\begin{lemma}\label{lemma:gpq} Let $\omega_1\in E_1^r,\omega_2\in E_1^s$ be effectives then
$$<g^p\omega_1,g^q\omega_2>=0\quad {\text if}\quad (p\not =q)
\quad {\rm or}
\quad ( p=q \quad {\rm and}\quad r\not=s)$$
Furthermore, in the case $p=q\geq 1$ and $r=s$, we have
$$<g^p\omega_1,g^p\omega_2>=p!(\prod_{i=0}^{p-1}(n-2r-i))<\omega_1,\omega_2>$$
\end{lemma}
{\bf Proof.} Recall that (see formula \ref{ck:glef})
$c^p(g^q\omega_2)=0$ if $p>q$, and  $c^q(g^p\omega_1)=0$ if $p<q$. This proves the
first part of the lemma. Also by the same formula and formula \ref{adj:gc}
 we have
$$<g^p\omega_1,g^p\omega_2>=<\omega_1,c^p(g^p\omega_2)>=<\omega_1,
p!(\prod_{i=0}^{p-1}(n-2r-i))\omega_2>$$
\ppp
\begin{corollary} Let $q=s+t$ then
$$h_{2q}={1\over (n-2q)!}\sum^{{\min\{2s,n-2s\}}}_{i=0}(-1)^i(n-2i)!<(R^s)_i,(R^t)_i>$$
In particular, we have
$$h_{4q}={1\over (n-4q)!}\sum^{{\min\{2q,n-2q\}}}_{i=0}(-1)^i
(n-2i)!<(R^q)_i,(R^q)_i>$$
\end{corollary}
{\bf Proof.} Note that
$$ h_{2q}={1\over (n-2q)!}*(g^{n-2q}R^q)=
 {1\over (n-2q)!}*(g^kR^sg^lR^t)$$
 where $k+l=n-2q$ and $s+t=q$. Then we apply
 the previous theorem to get
$$h_{2q}={1\over (n-2q)!}\sum^{{\min\{k+2s,l+2t\}}}_{i=0}(-1)^i
(n-2i)!<(g^kR^s)_i,(g^lR^t)_i>$$
Recall that $(g^kR^s)_i=R^s_i$ if $i\leq 2s$ otherwise it is zero. The same is true for
$(g^lR^t)_i$. This completes the proof.\ppp
\par\noindent\bigskip
The special case $q=1$ is of special interest. It provides an
 obstruction to the existence of an Einstein metric or a
  conformally flat metric with
zero scalar curvature in arbitrary  higher dimensions, as follows:  \par\medskip\noindent
\begin{theorem}\begin{enumerate}
\item If $(M,g)$ is an Einstein manifold with dimension $n\geq 4$,
then $h_4\geq 0$ and $h_4\equiv 0$ if and only if $(M,g)$ is flat.
\item If $(M,g)$ is a conformally flat  manifold with zero scalar curvature and
 dimension $n\geq 4$,
then $h_4\leq 0$ and $h_4\equiv 0$ if and only if $(M,g)$ is flat.
\end{enumerate}
\end{theorem}
{\bf Proof.} Straightforward using the previous corollary and theorem
\ref{sp+sp:sp-sp}.
\par\noindent\medskip
The previous theorem  can be generalized as follows. Its proof is also a
direct consequence of  the previous corollary and theorem \ref{pqcurv:thm}.
\begin{theorem} Let $(M,g)$ be a Riemannian manifold with dimension
 $n=2r\geq 4q$, for some $q\geq 1$.\begin{enumerate}
\item If $s_{(r,q)}(P)=s_{(r,q)}(P^\bot )$ for all $r$-planes $P$ then
 $h_{4q}\geq 0$
and $h_{4q}\equiv 0$  if and only if $(M,g)$ is flat.
\item If $s_{(r,q)}(P)=-s_{(r,q)}(P^\bot )$ for all $r$-planes $P$ then
 $h_{4q}\leq 0$ and $h_{4q}\equiv 0$ if and only if $(M,g)$ is flat.
 \end{enumerate}
 \end{theorem}
 The previous two theorems generalize similar results of Thorpe \cite{Thorpe}
 in the case $n=4q$.

\end{document}